\documentclass[12pt,a4paper]{article}
\usepackage{amsmath, amssymb, amsthm, amsfonts}

\numberwithin{equation}{section}

\begin{document}
\title{ On the $k$-free divisor problem}
\author{Jun FURUYA   and Wenguang ZHAI  }
\date{}

 \footnotetext[0]{The second-named author is
supported by National Natural Science Foundation of China(Grant
No. 10301018).} \maketitle

\begin{center}{Jun Furuya,\\
Department of Integrated Arts and Science,\\
Okinawa National College of Technology,\\
Nago, Okinawa, 905-2192,
Japan\\
E-mail: jfuruya@okinawa-ct.ac.jp}\end{center}

\begin{center}{
Wenguang Zhai,\\
School of Mathematical Sciences,\\
Shandong Normal University,\\
Jinan, Shandong, 250014,
P.R.China\\
E-mail: zhaiwg@hotmail.com}\end{center}

\begin{center}{Acta Arith.{\bf 123}(2006), 267-287}\end{center}

{\bf Abstract.}  Let $\Delta^{(k)}(x)$ denote the error term of the
$k$-free divisor problem for $k\geq 2$. In this paper   we establish
an asymptotic formula of the integral
$\int_1^T|\Delta^{(k)}(x)|^2dx$  for each $k\geq 4.$

{\section{Introduction}}

Let $d(n)$ denote the divisor function. Dirichlet first proved
that the error term
$$\Delta(x): =\sum_{n\leq x}d(n)-x\log x-
(2\gamma-1)x,\hspace{3mm}x\geq 2$$ satisfies
$\Delta(x)=O(x^{1/2}).$ The exponent $1/2$ was improved by many
authors. The latest result is due to Huxley [4], who proved that
\begin{equation*}
\Delta(x)= \left(x^{131/416}(\log x)^{26957/8320}\right).
\end{equation*}
It is conjectured that
\begin{equation}
\Delta(x)=O(x^{1/4+\varepsilon}),\label{1.2}
\end{equation}
which is supported by the classical mean-square result
\begin{equation}
\int_1^T\Delta^2(x)dx=\frac{(\zeta(3/2))^4}{6\pi^2
\zeta(3)}T^{3/2}+O(T\log^5 T)\label{1.3}
\end{equation}
proved by Tong [15].

Let $k\geq 2$ denote a fixed integer. An integer $n$ is called
$k$-free if $p^k$ does not divide $n$ for any prime $p$. Let
$d^{(k)}(n)$ denote the number of $k$-free divisors of the
positive integer $n$ and define
$$D^{(k)}(x):=\sum_{n\leq x}d^{(k)}(n).$$
Then the expected asymptotic formula of $D^{(k)}(x)$ is
\begin{equation*}
D^{(k)}(x)=C_1^{(k)}x\log x+C_2^{(k)}x+\Delta^{(k)}(x),
\end{equation*}
where $C_1^{(k)},C_2^{(k)}$ are two constants, $\Delta^{(k)}(x)$
is the error term. In 1874 Mertens [9] proved that
$\Delta^{(2)}(x)\ll x^{1/2}\log x.$ In 1932 H\"older \cite{s12}
proved that
\begin{eqnarray*}
\Delta^{(k)}(x)\ll\left\{\begin{array}{ll}
x^{1/2},&\mbox{if $k=2,$}\\
x^{1/3},& \mbox{if $k=3,$}\\
x^{33/100},&\mbox{if $k\geq 4.$}\\
\end{array}\right.
\end{eqnarray*}

For $k=2,3,$ it is very difficult to improve the exponent $1/k$ in
the bound $\Delta^{(k)}(x)\ll x^{1/k},$ unless we have substantial
progress in the study of the zero-free region of $\zeta(s).$
Therefore it is reasonable to get better improvements by assuming
the truth of the Riemann Hypothesis (RH). Such results were given
in [1, 2, 9, 12, 13, 14]. Especially in \cite{s2} R. C. Baker
proved $\Delta^{(2)}(x)\ll x^{4/11+\varepsilon}$ and in \cite{s8}
Kumchev proved $\Delta^{(3)}(x)\ll x^{27/85+\varepsilon}$ under
RH. For $k\geq 4,$ it is easy to show that if $\Delta(x)\ll
x^\alpha$ is true, then the estimate $\Delta^{(k)}(x)\ll
x^\alpha\log x$ follows.

We believe that the estimate
\begin{equation}
\Delta^{(k)}(x)=O( x^{1/4+\varepsilon})\label{1.5}
\end{equation}
would be true for any $k\geq 2$, which is an analogue of
(\ref{1.2}). For $k\geq 4$ it is easily seen that if the
conjecture (\ref{1.2}) is true, then so is (\ref{1.5}). For
$k=2,3,$ we cannot deduce the conjecture (\ref{1.5}) from
(\ref{1.2}) directly ; in this case we don't know the truth of
(\ref{1.5}) even if both (\ref{1.2}) and RH are true. However for
any $k\geq 2$, the conjecture (\ref{1.5}) cannot be proved by the
present method.

In this paper we shall study the mean square of $\Delta^{(k)}(x) $
for $k\geq 4,$ from which the truth of the conjecture (\ref{1.5})
$(k\geq 4)$ is supported partly. Our result is an analogue of
(\ref{1.3}).

{\bf Theorem 1.} We have the asymptotic formula
\begin{eqnarray*}
\int_1^T|\Delta^{(k)}(x)|^2dx=\frac{B_k}{6\pi^2}
T^{3/2}+\left\{\begin{array}{ll}
O(T^{3/2}e^{-c\delta(T)}),&\mbox{for $k=4$, }\\
O(T^{\delta_k+\varepsilon}),& \mbox{for $k\geq 5,$}
\end{array}\right.
\end{eqnarray*}
where
\begin{align*}
&B_k:=\sum_{m=1}^\infty g_k^2(m)m^{-3/2},\quad
g_k(m): =\sum_{m=nd^k}\mu(d)d(n)d^{k/2}, \\
&\delta(u): =(\log u)^{3/5}(\log\log u)^{-1/5},\\
&\delta_5:=75/52,
\hspace{2mm}\delta_k:=3/2-1/2k+1/k^2\hspace{2mm}(k\geq 6),
\end{align*}
 and where $c>0$ is an absolute constant .

{\bf Corollary 1.} If $k\geq 4,$ then we have
$$\Delta^{(k)}(x)=\Omega(x^{1/4}).$$

By the same method we can study the mean square of
$\Delta(1,1,k;x),$ which is defined by
$$\Delta(1,1,k;x):=\sum_{n\leq x}d(1,1,k;n)-x\left\{\zeta(k)\log
x+k\zeta^{\prime}(k)+(2\gamma-1)\zeta(k)\right\}
-\zeta^2\big(\frac1k\big)x^{1/k},$$ where
$d(1,1,k;n)=\sum_{n=m_1m_2d^k}1$ and $\gamma$ is the Euler
constant. This is a special three-dimensional divisor problem.
From the formula (5.3) of Ivi\'c \cite{s6} we have
\begin{equation}
\int_1^T\Delta^2(1,1,k;x)dx\ll T^{3/2+\varepsilon}.
\end{equation}
From Kr\"atzel \cite{s7} we know that
\begin{equation}\Delta(1,1,k;x)=\Omega(x^{1/4})\end{equation} if $k\geq 5.$

Now we  prove the following Theorem 2, which improves (1.4).

{\bf Theorem 2.} Suppose $k\geq 3$ is a fixed integer. Then we
have
\begin{eqnarray*}
\int_1^T\Delta^2(1,1,k;x)dx=\frac{C_k}{6\pi^2}
T^{3/2}+\left\{\begin{array}{ll}
O(T^{53/36}\log^3 T),&\mbox{if $k=3$, }\\
O(T^{29/20}\log^{503}T),& \mbox{if $k=4,$}\\
O(T^{75/52}\log^{1000}T),& \mbox{if $k=5,$}\\
O(T^{3/2-1/2k+1/k^2+\varepsilon}),& \mbox{if $k\geq 6,$}\\
\end{array}\right.
\end{eqnarray*}
where
\begin{eqnarray*}
&&C_k: =\sum_{m=1}^\infty f_k^2(m)m^{-3/2},\quad f_k(m):
=\sum_{m=nd^k}d(n)d^{k/2}.
\end{eqnarray*}

{\bf Corollary 2.} The formula (1.5) holds for $k=3,4.$

{\bf Notations.} For a real number $u,$ $[u]$ denotes the integer
part of $u,$ $\{u\}$ denotes the fractional part of $u,$
$\psi(u)=\{u\}-1/2,$ $\Vert u\Vert$ denotes the distance from $u$
to the integer nearest to $u.$ $\mu(d)$ is the M\"obius function.
Let $(m,n)$ denote the greatest common divisor of natural numbers
$m$ and $n.$ $n\sim N$ means $N<n\leq 2N.$ $\varepsilon$ always
denotes a sufficiently small positive constant which may be
different at different places. $SC(\Sigma)$ denotes the summation
condition of the sum $\Sigma.$

{\section{The expression of $\Delta^{(k)}(x)$}}

In order to prove Theorem 1, we shall give a simple expression of
$\Delta^{(k)}(x)$ in this section.

{\bf Lemma 2.1.} There exists an absolute constant $c_1>0$ such
that the estimate
$$M(u): =\sum_{n\leq u}\mu(n)\ll ue^{-c_1\delta(u)}$$
holds for $u\geq 2.$

This is Theorem 12.7 of Ivi\'c \cite{s5}. Now we prove the
following

{\bf Lemma 2.2.} Suppose $10\leq y\ll x^{1/k},$ then we have
$$\Delta^{(k)}(x)=\sum_{d\leq y}\mu(d)\Delta\big(\frac{x}{d^k}\big)
+O\left(xy^{1-k}e^{-c_1\delta(y)}\log x\right).$$
\begin{proof}
We have
\begin{align*}
D^{(k)}(x)=&\ \sum_{\stackrel{mn\leq x}{\scriptstyle
m:k{\mbox{-}}free}}1
=\sum_{d^kmn\leq x}\mu(d)=\sum_{d^kn\leq x}\mu(d)d(n)\\
=&\ \sum_{d\leq y}\mu(d)D\big(\frac{x}{d^k}\big)+\sum_{n\leq
x/y^k}d(n)
M\left(\big(\frac{x}{n}\big)^{1/k}\right)-D\big(\frac{x}{y^k}\big)M(y)\\
=&\ {\sum}_1+{\sum}_2-{\sum}_3,
\end{align*}
say. From Lemma 2.1 and the estimate $D(u)\ll u\log u$ directly we
have
$${\sum}_3\ll xy^{1-k}e^{-c_1\delta(y)}\log x.$$
From Lemma 2.1, the estimate $D(u)\ll u\log u$ and partial
summation we have
$${\sum}_2\ll xy^{1-k}e^{-c_1\delta(y)}\log x$$
if we note that $e^{-c_1\delta((x/n)^{1/k})}\leq
e^{-c_1\delta(y)}$ for all $n\leq x/y^k.$ By Lemma 2.1 and simple
calculations we have
\begin{align*}
{\sum}_1=&\ \sum_{d\leq
y}\mu(d)\left\{\frac{x}{d^k}\log\frac{x}{d^k}+(2\gamma-1)
\frac{x}{d^k}\right\}
+\sum_{d\leq y}\mu(d)\Delta\big(\frac{x}{d^k}\big)\\
=&\ {\mbox{(Main term)}}+\sum_{d\leq
y}\mu(d)\Delta\big(\frac{x}{d^k}\big)
+O\left(xy^{1-k}e^{-c_1\delta(y)}\log x\right).
\end{align*}
Whence Lemma 2.2 follows.\end{proof}

{\section{Proof of Theorem 1(Beginning)}}

Suppose $T\geq 10$ is large. It suffices for us to evaluate the
integral $\int_T^{2T}|\Delta^{(k)}(x)|^2dx.$

Let $T^\varepsilon\ll y\ll T^{1/k-\varepsilon},$ $T^\varepsilon\ll
z\ll T^{1-\varepsilon}$ be two parameters to be determined later.
Let
$$\Delta_1(u):=\frac{u^{1/4}}{\pi\sqrt 2}\sum_{n\leq z}\frac{d(n)}{n^{3/4}}
\cos\left(4\pi\sqrt{nu}-\frac{\pi}{4}\right),\
\Delta_2(u;z):=\Delta(u)-\Delta_1(u).$$ Then by Lemma 2.2 we can
write
\begin{equation}
\Delta^{(k)}(x)=R_1^{(k)}(x)+R_2^{(k)}(x)+R_3^{(k)}(x),
\label{3.1}
\end{equation}
where
\begin{align*}
R_1^{(k)}(x):=&\ \frac{x^{1/4}}{\pi\sqrt 2}\sum_{d\leq
y}\frac{\mu(d)}{d^{k/4}}\sum_{n\leq z}\frac{d(n)}{n^{3/4}}
\cos\left(4\pi\sqrt{\frac{nx}{d^k}}-\frac{\pi}{4}\right),\\
R_2^{(k)}(x):=&\ \sum_{d\leq
y}\mu(d)\Delta_2(\frac{x}{d^k};z)\quad {\mbox{and}}\quad
R_3^{(k)}(x):=O\left(xy^{1-k}e^{-c_1\delta(y)}\log x\right).
\end{align*}

{\bf Lemma 3.1.} Suppose $A>0$ is any fixed constant,
$T^\varepsilon\ll V\ll T^A.$ Then we have
\begin{eqnarray*}
\int_V^{2V}\Delta_2^2(u;z)du\ll V^{3/2}z^{-1/2}\log^3 V+ V\log^5
V.
\end{eqnarray*}

\begin{proof} Suppose $\min(z,V^{11})<N\ll V^{B} $ is a large
parameter, where $B>0$ is a constant suitably large. By Lemma 3 of
Meurman \cite{s10} we have
\begin{equation*}
\Delta_2(u;z)=\frac{u^{1/4}}{\pi\sqrt 2}\sum_{z<n\leq
N}\frac{d(n)}{n^{3/4}}
\cos\left(4\pi\sqrt{nu}-\frac{\pi}4\right)+\Delta_2(u;N),
\end{equation*}
where $\Delta_2(u;N)\ll u^{-1/4}$ if $\Vert u\Vert\gg
u^{5/2}N^{-1/2},$ and $\Delta_2(u;N)\ll u^{\varepsilon}$
otherwise. Thus we have
\begin{eqnarray*}
&&\int_V^{2V}\Delta_2^2(u;z)du\ll \int_1+\int_2,
\end{eqnarray*}
where
\begin{eqnarray*}
&&\int_1=\int_V^{2V}\left|u^{1/4}\sum_{z<n\leq
N}\frac{d(n)}{n^{3/4}}\cos\left(4\pi\sqrt{nu}-\frac{\pi}4\right)\right|^2du,
\quad\int_2=\int_V^{2V}\Delta_2^2(u;N)du
\end{eqnarray*}
For $\int_1$ we have
\begin{eqnarray*}
\int_1 &&\ll\int_V^{2V}|u^{\frac{1}{4}}\sum_{z<n\leq
N}\frac{d(n)}{n^{ 3/4}} e(2\sqrt{nu})|^2du\\
&&\ll T^{3/2}\sum_{z<n\leq N}\frac{d^2(n)}{n^{3/2}}
+V\sum_{z<m<n\leq N}\frac{d(n)d(m)}{(mn)^{3/4}(\sqrt n-\sqrt m)}\\
&& \ll \frac{V^{3/2}\log^3 V}{z^{1/2}}+V\log^5 V,
\end{eqnarray*}
where we used the well-known estimates
\begin{align}
&\sum_{n\leq u}d^2(n)\ll u\log^3 u,\label{3.3}\\
&\sum_{z<m<n\leq N}\frac{d(n)d(m)}{(mn)^{3/4}(\sqrt n-\sqrt m)}\ll
\log^5 N\ll \log^5 V.\nonumber
\end{align}
For $\int_2$ we have
\begin{eqnarray*}
\int_2\ll V(V^{5/2+\varepsilon}N^{-1/2}+V^{-1/4})\ll
V^{7/2+\varepsilon}N^{-1/2}+V^{3/4}\ll V.
\end{eqnarray*}
Now Lemma 3.1 follows from the above estimates.
\end{proof}

By Cauchy's inequality and Lemma 3.1 we get
\begin{eqnarray}
&&\int_T^{2T}|R_2^{(k)}(x)|^2dx=\int_T^{2T}\left|\sum_{d\leq
y}\mu(d)d^{-1/2}d^{1/2}\Delta_2\left(\frac{x}{d^k};z\right)\right|^2dx\\
&&\ll \int_T^{2T}\left(\sum_{d\leq
y}d^{-1}\right)\left(\sum_{d\leq
y}d\left|\Delta_2\left(\frac{x}{d^k};z\right)\right|^2\right)dx\nonumber\\
&&\ll \sum_{d\leq
y}d\int_T^{2T}\left|\Delta_2\left(\frac{x}{d^k};z\right)\right|^2dx\log
y\nonumber\\
&&\ll \sum_{d\leq
y}d^{k+1}\int_{T/d^k}^{2T/d^k}\left|\Delta_2(u;z)\right|^2du\log
y\nonumber\\
&&\ll \sum_{d\leq
y}d^{k+1}\left(\left(\frac{T}{d^k}\right)^{3/2}z^{-1/2}\log^3T+Td^{-k}\log^5T\right)\log
y\nonumber\\
&&\ll T^{3/2}z^{-1/2}\sum_{d\leq y}d^{1-k/2}\log^4 T +Ty^2\log^6
T\nonumber\\
&&\ll\left\{\begin{array}{ll}
T^{3/2}z^{-1/2}y^{1/2}\log^4 T +Ty^2\log^6 T,&\mbox{if $k=3$, }\\
T^{3/2}z^{-1/2}\log^5 T +Ty^2\log^6T ,& \mbox{if
$k\geq4.$}\nonumber
\end{array}\right.
\end{eqnarray}

If $k=4,$ we take $y=T^{1/4}e^{-c_2\delta(T)},$ where
$c_2=c_1/4^{8/5}.$  It is easy to see that $R_3^{(k)}(x)\ll
T^{1/4}e^{-c_3\delta(T)}$ holds for all $T\leq x\leq 2T$, where
$0<c_3<c_1/4^{8/5}$ is an absolute constant. Hence
\begin{equation}
\int_T^{2T}|R_3^{(4)}(x)|^2dx\ll T^{3/2}e^{-2c_3\delta(T)}.
\end{equation}

If $k\geq 5,$ then  we have
\begin{equation}
\int_T^{2T}|R_3^{(k)}(x)|^2dx\ll T^{3}y^{2-2k}.
\end{equation}

Now we consider the mean square of $R_1^{(k)}(x).$ By the
elementary formula
$$\cos u\cos v=\frac{1}{2}(\cos {(u-v)}+\cos {(u+v)})$$
we may write
\begin{align}
|R_1^{(k)}(x)|^2&=\frac{x^{1/2}}{2\pi^2}\sum_{d_1,d_2\leq y}
\frac{\mu(d_1)\mu(d_2)}{(d_1d_2)^{k/4}}\sum_{n_1,n_2\leq z}
\frac{d(n_1)d(n_2)}{(n_1n_2)^{3/4}}\\
&\hspace{10mm}\times
\cos{\left(4\pi\sqrt{\frac{n_1x}{d_1^k}}-\frac{\pi}{4}\right)}
\cos{\left(4\pi\sqrt{\frac{n_2x}{d_2^k}}-\frac{\pi}{4}\right)}\nonumber\\
&=S_1(x)+S_2(x)+S_3(x),\nonumber
\end{align}
where
\begin{align*}
S_1(x)=&\ \frac{x^{1/2}}{4\pi^2}\sum_{\stackrel{d_1,d_2\leq
y;n_1,n_2\leq z}{n_1d_2^k=n_2d_1^k}}
\frac{\mu(d_1)\mu(d_2)}{(d_1d_2)^{k/4}}
\frac{d(n_1)d(n_2)}{(n_1n_2)^{3/4}},\\
S_2(x)=&\ \frac{x^{1/2}}{4\pi^2}\sum_{\stackrel{d_1,d_2\leq
y;n_1,n_2\leq z}{n_1d_2^k\not =n_2d_1^k}}
\frac{\mu(d_1)\mu(d_2)}{(d_1d_2)^{k/4}}
\frac{d(n_1)d(n_2)}{(n_1n_2)^{3/4}}\cos{\left(4\pi\sqrt x
\left(\sqrt
{\frac{n_1}{d_1^k}}-\sqrt{\frac{n_2}{d_2^k}} \right) \right)},\\
S_3(x)=&\ \frac{x^{1/2}}{4\pi^2}\sum_{d_1,d_2\leq y;n_1,n_2\leq z}
\frac{\mu(d_1)\mu(d_2)}{(d_1d_2)^{k/4}}
\frac{d(n_1)d(n_2)}{(n_1n_2)^{3/4}}\sin{\left(4\pi\sqrt x
\left(\sqrt {\frac{n_1}{d_1^k}}+\sqrt{\frac{n_2}{d_2^k}} \right)
\right)}.
\end{align*}

We have
\begin{align}
&\int_T^{2T}S_1(x)dx=\frac{B_k(y,z)}{4\pi^2}\int_T^{2T}x^{1/2}dx,\\
& B_k(y,z):=\sum_{\stackrel{d_1,d_2\leq y;n_1,n_2\leq
z}{n_1d_2^k=n_2d_1^k}} \frac{\mu(d_1)\mu(d_2)}{(d_1d_2)^{k/4}}
\frac{d(n_1)d(n_2)}{(n_1n_2)^{3/4}}.\nonumber
\end{align}

By the first derivative test we get
\begin{equation}
\int_T^{2T}S_2(x)dx\ll TE_k(y,z),
\end{equation}
where
$$E_k(y,z)=\sum_{\stackrel{d_1,d_2\leq
y;n_1,n_2\leq z}{n_1d_2^k\not =n_2d_1^k}}
\frac{d(n_1)d(n_2)}{(d_1d_2)^{k/4}(n_1n_2)^{3/4}}
\min\left(T^{1/2},\frac{1}{\left|\sqrt
{\frac{n_1}{d_1^k}}-\sqrt{\frac{n_2}{d_2^k}}\right|}\right).$$

By the first derivative test again we get
\begin{align}
\int_T^{2T}S_3(x)dx&\ll \sum_{d_1,d_2\leq y;n_1,n_2\leq z}
\frac{d(n_1)d(n_2)}{(d_1d_2)^{k/4}(n_1n_2)^{3/4}}\frac{1}{\left|\sqrt
{\frac{n_1}{d_1^k}}+\sqrt{\frac{n_2}{d_2^k}}\right|}\\
&\ll \sum_{d_1,d_2\leq y;n_1,n_2\leq z}
\frac{d(n_1)d(n_2)}{(d_1d_2)^{k/4}(n_1n_2)^{3/4}}\frac{1}{\left(\sqrt
{\frac{n_1}{d_1^k}}\sqrt {\frac{n_2}{d_2^k}}\right)^{1/2}}\nonumber\\
&\ll \sum_{d_1,d_2\leq y;n_1,n_2\leq z}
\frac{d(n_1)d(n_2)}{n_1n_2}\ll y^2\log^4z,\nonumber
\end{align}
where the inequality $ab\geq 2\sqrt{ab}$ and the estimate $D(u)\ll
u\log u$ were used.

Now the problem is reduced to evaluating $B_k(y,z)$ and estimating
$E_k(y,z).$

\section{\bf Evaluation of $B_k(y,z)$}

In this section we shall evaluate $B_k(y,z)$. We have
\begin{align*}
B_k(y,z)&=\sum_{\stackrel{d_1,d_2\leq y;n_1,n_2\leq
z}{n_1d_2^k=n_2d_1^k}}
\frac{\mu(d_1)\mu(d_2)d(n_1)d(n_2)(d_1d_2)^{k/2}}{(n_1d_2^kn_2d_1^k)^{3/4}}\\
&=\sum_{m\leq zy^k}g^2(m;y,z)m^{-3/2},\nonumber
\end{align*}
where
$$g_k(m;y,z):=\sum_{\stackrel{m=nd^k}{n\leq z,d\leq
y}}\mu(d)d(n)d^{k/2}.$$ Let
$$g_k(m)=\sum_{m=nd^k}\mu(d)d(n)d^{k/2},\hspace{2mm}g_0(m)=f_k(m)=\sum_{m=nd^k}d(n)d^{k/2}.$$
Let $z_0: =\min(y,z).$ Obviously,
\begin{eqnarray*}
&&g_k(m;y,z)=g_k(m),\hspace{2mm}m\leq z_0,\\
&&|g_k(m;y,z)|\leq g_0(m), |g_k(m)|\leq g_0(m), m\geq 1.
\end{eqnarray*}

Thus
\begin{align}
B_k(y,z)&=\sum_{m\leq z_0}g_k^2(m)m^{-3/2}+\sum_{z_0<m\leq
zy^k}g_k^2(m;y,z)m^{-3/2}\label{4.2}\\
&=\sum_{m\leq z_0}g_k^2(m)m^{-3/2}+O\left(\sum_{z_0<m\leq
zy^k}|g_0^2(m)|m^{-3/2}\right)\nonumber.
\end{align}

For any $1<U<V<\infty,$ we shall estimate the sum
$$W_k(U,V):= \sum_{U<m\leq
V}|g_0^2(m)|m^{-3/2}.$$

Obviously $g_0(m)$ is a multiplicative function. So for $m> 1,$ we
have
$$g_0(m)=\prod_{p^\alpha \Vert m}g_0(p^\alpha).$$
If $1\leq \alpha\leq k-1,$ then $g_0(p^\alpha)=\alpha+1,$ which
implies that if $n$ is $k$-free then $g_0(n)=d(n).$

Now suppose $ek\leq \alpha<(e+1)k$ for some integer $e\geq 1.$ It
can be easily seen that if we write $p^\alpha$ in the form
$p^\alpha=nd^k,$ then $n=p^{\alpha-jk},d=p^{j}, j=0,1,2,\cdots,e.$
Then we have
\begin{align*}
g_0(p^\alpha)=&\
\sum_{j=0}^{e}(\alpha-jk+1)p^{jk/2}=p^{ek/2}\sum_{j=0}^{e}(\alpha-jk+1)p^{-(e-j)k/2}\\
\leq&\
(\alpha+1)p^{ek/2}\sum_{j=0}^{e}p^{-(e-j)k/2}=(\alpha+1)p^{ek/2}\sum_{j=0}^{e}p^{-jk/2}\\
\leq&\ (\alpha+1)p^{ek/2}\sum_{j=0}^{\infty}2^{-jk/2}\leq
2(\alpha+1)p^{\alpha/2},
\end{align*}
which implies that if $l$ is $k$-full, then
$$g_0(l)\leq \prod_{p^\alpha \Vert l}2(\alpha+1)p^{\alpha/2}
=2^{\omega(l)}d(l)l^{1/2}\leq d^2(l)l^{1/2}.$$

Let $\delta_{(k)}(n),$ $\delta^{(k)}(n)$ denote the characteristic
function of $k$-free and $k$-full numbers, respectively. Each
integer $m$ can be uniquely written as $m=nl,$ $(n,l)=1,$
$\delta_{(k)}(n)=1,\delta^{(k)}(l)=1.$ Thus we have
\begin{align*}
W_k(U,V)=&\ \sum_{\stackrel{U<nl\leq
V}{(n,l)=1}}g_0^2(n)g_0^2(l)\delta_{(k)}(n)\delta^{(k)}(l)(nl)^{-3/2}\\
\ll&\ {\sum}_4+{\sum}_5,
\end{align*}
where
\begin{eqnarray*}
&&{\sum}_4:=\sum_{l\leq U/3, U<nl\leq
V}g_0^2(n)g_0^2(l)\delta_{(k)}(n)\delta^{(k)}(l)(nl)^{-3/2},\\
&&{\sum}_5:=\sum_{l> U/3, U<nl\leq
V}g_0^2(n)g_0^2(l)\delta_{(k)}(n)\delta^{(k)}(l)(nl)^{-3/2}.
\end{eqnarray*}

{\bf Lemma 4.1.} We have the estimate
\begin{eqnarray}
&&\sum_{n\leq u}d^4(n)\delta^{(k)}(n)\ll u^{1/k}\log^{(k+1)^4-1}
u, u\geq 2.\label{4.3}
\end{eqnarray}
\begin{proof}
For $\Re s>1/k,$ it is easy to show that
$$\sum_{n=1}^\infty
d^4(n)\delta^{(k)}(n)n^{-s}=\zeta^{(k+1)^4}(ks)G_k(s),$$ where
$G_k(s)$ is absolutely convergent for $\Re s>1/(1+k).$ And whence
(\ref{4.3}) follows.
\end{proof}

By (\ref{3.3}), partial summation and Lemma 4.1 we have
\begin{align*}
{\sum}_4\ll&\ \sum_{l\leq
U/3}g_0^2(l)\delta^{(k)}(l)l^{-3/2}\sum_{U/l<n\leq
V/l}g_0^2(n)n^{-3/2}\\
\ll&\ \sum_{l\leq
U/3}g_0^2(l)\delta^{(k)}(l)l^{-3/2}(U/l)^{-1/2}\log^3 U\\
\ll&\ U^{-1/2}\log^3 U\sum_{l\leq U/3}d^4(l)\delta^{(k)}(l)\\
\ll&\ U^{-1/2+1/k}\log^{(k+1)^4+2} U
\end{align*}
and
\begin{align*}
{\sum}_5\ll&\ \sum_{l>
U/3}g_0^2(l)\delta^{(k)}(l)l^{-3/2}\sum_{m}g_0^2(n)n^{-3/2}\\
\ll&\ \sum_{l>
U/3}d^4(l)\delta^{(k)}(l)l^{-1/2}\\
\ll&\ U^{-1/2+1/k}\log^{(k+1)^4+2} U.
\end{align*}
Thus
\begin{equation}
W_k(U,V)\ll U^{-1/2+1/k}\log^{(k+1)^4+2} U.\label{4.4}
\end{equation}

From (\ref{4.2}) and (\ref{4.4}) we immediately get
\begin{equation}
B_k(y,z)=\sum_{m=1}^\infty
g_k^2(m)m^{-3/2}+O\left(z_0^{-1/2+1/k}\log^{(k+1)^4+2} z_0\right).
\end{equation}

\section{\bf Estimation of $E_k(y,z)$}

In this section we shall estimate $E_k(y,z).$ By a splitting
argument, we have
\begin{equation}
E_k(y,z)\ll E_k(D_1,D_2,N_1,N_2)z^\varepsilon \log^2 y
\end{equation}
for some $(D_1,D_2,N_1,N_2)$ with $1\ll D_j\ll y,1\ll N_j\ll z,
j=1,2,$ where
\begin{eqnarray*}
&&E_k(D_1,D_2,N_1,N_2)=\sum
\frac{1}{(d_1d_1)^{k/4}(n_1n_2)^{3/4}}\min\left(T^{1/2},\frac{1}{|\sqrt
{\frac{n_1}{d_1^k}}-\sqrt{\frac{n_2}{d_2^k}}|}\right),\\
&&SC(\sum): d_1\sim D_1, d_2\sim D_2,n_1\sim N_1, n_2\sim N_2,
n_1d_2^k\not= n_2d_1^k.
\end{eqnarray*}
We write
\begin{align*}
E_k(D_1,D_2,N_1,N_2)=&\ {\sum}_6
\frac{1}{(d_1d_1)^{k/4}(n_1n_2)^{3/4}}\min\left(T^{1/2},\frac{1}{|\sqrt
{\frac{n_1}{d_1^k}}-\sqrt{\frac{n_2}{d_2^k}}|}\right)\\
&+{\sum}_7
\frac{1}{(d_1d_1)^{k/4}(n_1n_2)^{3/4}}\min\left(T^{1/2},\frac{1}{|\sqrt
{\frac{n_1}{d_1^k}}-\sqrt{\frac{n_2}{d_2^k}}|}\right),
\end{align*}
where
\begin{eqnarray*}
&&SC({\sum}_6): d_1\sim D_1, d_2\sim D_2,n_1\sim N_1, n_2\sim
N_2,\\&&\hspace{20mm} \left|\sqrt
{\frac{n_1}{d_1^k}}-\sqrt{\frac{n_2}{d_2^k}}\right|\geq
\left(\sqrt
{\frac{n_1}{d_1^k}}\sqrt{\frac{n_2}{d_2^k}}\right)^{1/2}/10,\\
&&SC({\sum}_7): d_1\sim D_1, d_2\sim D_2,n_1\sim N_1, n_2\sim
N_2,\\&&\hspace{20mm} \left|\sqrt
{\frac{n_1}{d_1^k}}-\sqrt{\frac{n_2}{d_2^k}}\right|< \left(\sqrt
{\frac{n_1}{d_1^k}}\sqrt{\frac{n_2}{d_2^k}}\right)^{1/2}/10.
\end{eqnarray*}
Trivially we have
\begin{equation}
{\sum}_6\ll \sum_{\stackrel{d_j\sim D_j,n_j\sim
N_j}{j=1,2}}\frac{1}{(d_1d_1)^{k/4}(n_1n_2)^{3/4}}\left(\sqrt
{\frac{n_1}{d_1^k}}\sqrt{\frac{n_2}{d_2^k}}\right)^{-1/2}\ll
D_1D_2\ll y^2.\label{5.2}
\end{equation}

Suppose $\delta>0,$ and let ${\cal A}(D_1,D_2,N_1, N_2;\delta)$
denote the number of the solutions of the inequality
\begin{equation}
\left|\sqrt
{\frac{n_1}{d_1^k}}-\sqrt{\frac{n_2}{d_2^k}}\right|\leq
\delta,d_1\sim D_1, d_2\sim D_2,n_1\sim N_1, n_2\sim N_2.
\label{5.3}
\end{equation}
In order to estimate $\sum_7,$ we need an upper bound of ${\cal
A}(D_1,D_2,N_1, N_2;\delta).$

{\bf Lemma 5.1.} We have
\begin{align*}
{\cal A}(D_1,D_2,N_1, N_2;\delta) &\ll\delta
(D_1D_2)^{1+k/4}(N_1N_2)^{3/4}\\
&+(D_1D_2N_1N_2)^{1/2}\log 2D_1D_2N_1N_2,
\end{align*} where
the implied constant is absolute.

\begin{proof}
We shall use an idea of Fouvry and Iwaniec [3]. Suppose $u$ and
$v$ are two positive integers and let ${\cal A}_{u,v}(D_1,D_2,N_1,
N_2;\delta)$ denote the number of the solutions of the inequality
(\ref{5.3}) with $(n_1,n_2)=u,(d_1,d_2)=v.$ Set $n_j=m_ju,d_j=l_jv
(j=1,2),$ then $(m_1,m_2,l_1,l_2)$ satisfies
\begin{equation}
\left|\sqrt{\frac{m_1}{m_2}}-\sqrt{\frac{l_1^k}{l_2^k}}\right|\leq
2^{k/2}\delta D_1^{k/2}N_2^{-1/2}
\end{equation}
and
\begin{equation}
\left|\sqrt{\frac{m_2}{m_1}}-\sqrt{\frac{l_2^k}{l_1^k}}\right|\leq
2^{k/2}\delta D_2^{k/2}N_1^{-1/2}.
\end{equation}

It is easy to show that $\sqrt{\frac{m_1}{m_2}}$ is
$u^2N_2^{-3/2}N_1^{-1/2}-$spaced, so from (5.4) we get
\begin{align*}
{\cal A}_{u,v}(D_1,D_2,N_1, N_2;\delta)&\ll
\frac{D_1D_2}{v^2}\large(1+\frac{\delta
D_1^{k/2}N_2N_1^{1/2}}{u^2}\large)\\
&\ll \frac{D_1D_2}{v^2}+\frac{\delta
D_1D_2D_1^{k/2}N_2N_1^{1/2}}{u^2v^2}.
\end{align*}
Similarly, since $\sqrt{\frac{m_2}{m_1}}$ is
$u^2N_1^{-3/2}N_2^{-1/2}-$spaced , from (5.5) we get
\begin{eqnarray*}
{\cal A}_{u,v}(D_1,D_2,N_1, N_2;\delta) \ll
\frac{D_1D_2}{v^2}+\frac{\delta
D_1D_2D_2^{k/2}N_1N_2^{1/2}}{u^2v^2}.
\end{eqnarray*}
From the above two estimates we get
\begin{align}
{\cal A}_{u,v}(D_1,D_2,N_1, N_2;\delta)& \ll
\frac{D_1D_2}{v^2}+\frac{\delta
D_1D_2}{u^2v^2}\min(D_1^{k/2}N_2N_1^{1/2},D_2^{k/2}N_1N_2^{1/2})
\nonumber\\
&\ll \frac{D_1D_2}{v^2}+\frac{\delta
(D_1D_2)^{1+k/4}(N_1N_2)^{3/4}}{u^2v^2}
\end{align}
if we note that $\min(a,b)\leq a^{1/2}b^{1/2}.$

It is easy to show that $(l_1/l_2)^{k/2}$ is
$v^2D_2^{-2}(D_1/D_2)^{k/2-1}-$spaced, from (5.4) we get
\begin{align*}
{\cal A}_{u,v}(D_1,D_2,N_1, N_2;\delta) & \ll
\frac{N_1N_2}{u^2}\large(1+\delta
D_1^{k/2}N_2^{-1/2}v^{-2}D_2^{2}(D_1/D_2)^{-k/2+1}\large)\\
&\ll \frac{N_1N_2}{u^2}+\frac{\delta
D_1D_2D_2^{k/2}N_1N_2^{1/2}}{u^2v^2}.
\end{align*}
Similarly from (5.5) we get
\begin{eqnarray*}
{\cal A}_{u,v}(D_1,D_2,N_1, N_2;\delta) \ll
\frac{N_1N_2}{u^2}+\frac{\delta
D_1D_2D_1^{k/2}N_2N_1^{1/2}}{u^2v^2}.\end{eqnarray*}

From the above two estimates we have
\begin{eqnarray*}
{\cal A}_{u,v}(D_1,D_2,N_1, N_2;\delta) \ll
\frac{N_1N_2}{u^2}+\frac{\delta
(D_1D_2)^{1+k/4}(N_1N_2)^{3/4}}{u^2v^2},\end{eqnarray*} which
combining (5.6) gives
\begin{eqnarray*}
{\cal A}_{u,v}(D_1,D_2,N_1, N_2;\delta) \ll \frac{\delta
(D_1D_2)^{1+k/4}(N_1N_2)^{3/4}}{u^2v^2}+
\min(\frac{N_1N_2}{u^2},\frac{D_1D_2}{v^2}).
\end{eqnarray*}
Summing over $u$ and $v$ completes the proof of Lemma 5.1.
\end{proof}

Now we estimate $\sum_7.$ Let $\Omega=\sqrt
{\frac{n_1}{d_1^k}}-\sqrt{\frac{n_2}{d_2^k}}.$ By Lemma 5.1 the
contribution of $T^{1/2}$ is (note that $|\Omega|\leq T^{-1/2}$)
\begin{eqnarray*}
&&\ll \frac{T^{1/2}}{(D_1D_2)^{k/4}(N_1N_2)^{3/4}} {\cal
A}_{u,v}(D_1,D_2,N_1, N_2;T^{-1/2})\\
&&\ll \frac{T^{1/2}\log T}{(D_1D_2)^{k/4-1/2}(N_1N_2)^{1/4}}
+D_1D_2.\end{eqnarray*}

Divide the remaining range into $O(\log T)$ intervals of the form
$T^{-1/2}<\delta<|\Omega|\leq 2\delta.$ By Lemma 5.1 again we find
that the contribution of $1/|\Omega|$ is
\begin{eqnarray*}
&&\ll \log T \max_{\delta>T^{-1/2}}\frac{{\cal
A}_{u,v}(D_1,D_2,N_1,
N_2;2\delta)}{(D_1D_2)^{k/4}(N_1N_2)^{3/4}\delta}
\\
&&\ll \frac{T^{1/2}\log^2 T}{(D_1D_2)^{k/4-1/2}(N_1N_2)^{1/4}}
+D_1D_2\log T.
\end{eqnarray*}

From the above two estimates we get
\begin{equation}
{\sum}_7\ll \frac{T^{1/2}\log^2
T}{(D_1D_2)^{k/4-1/2}(N_1N_2)^{1/4}} +y^2\log T.
\end{equation}

Now we give another estimate of $\sum_7.$ By noting that $\sqrt
{\frac{n_1}{d_1^k}}\asymp \sqrt{\frac{n_2}{d_2^k}}$ we get
\begin{align*}
\frac{1}{|\Omega|}=&\ \frac{\sqrt {\frac{n_1}{d_1^k}}+
\sqrt{\frac{n_2}{d_2^k}}}{|\frac{n_1}{d_1^k}- \frac{n_2}{d_2^k}|}
\ll \frac{(d_1d_2)^k\large(\sqrt {\frac{n_1}{d_1^k}}+
\sqrt{\frac{n_2}{d_2^k}}\large)}{|n_1d_2^k- n_2d_1^k|} \\
\ll&\ (d_1d_2)^k (\sqrt {\frac{n_1}{d_1^k}}\sqrt
{\frac{n_2}{d_2^k}})^{1/2} \ll (d_1d_2)^{3k/4}(n_1n_2)^{1/4}\\
\ll&\ (D_1D_2)^{3k/4}(N_1N_2)^{1/4}.
\end{align*}
The range of $\Omega$ can be divided into $O(\log T)$ intervals of
the form $$(D_1D_2)^{-3k/4}(N_1N_2)^{-1/4}\ll \delta\leq
|\Omega|\leq 2\delta.$$ By Lemma 5.1 we have
\begin{align}
{\sum}_7&\ll \frac{1}{(D_1D_2)^{k/4}(N_1N_2)^{3/4}}
\sum_{\Omega}\frac{1}{|\Omega|}\\
&\ll \frac{\log
T}{(D_1D_2)^{k/4}(N_1N_2)^{3/4}}\max_{\delta}\frac{{\cal
A}_{u,v}(D_1,D_2,N_1, N_2;\delta)}{\delta}\nonumber\\
&\ll (D_1D_2)^{(k+1)/2}\log^2 T\nonumber
\end{align}
if we note that $\delta\gg (D_1D_2)^{-3k/4}(N_1N_2)^{-1/4}.$

From (5.7) and (5.8) we get
\begin{align}
\sum_7&\ll y^2\log
T+\min\left(\frac{T^{1/2}}{(D_1D_2)^{k/4-1/2}(N_1N_2)^{1/4}}
,(D_1D_2)^{(k+1)/2}\right)\log^2 T\\
&\ll y^2\log T+
\left(\frac{T^{1/2}}{(D_1D_2)^{k/4-1/2}(N_1N_2)^{1/4}}\right)^{(2k+2)/3k}
\left((D_1D_2)^{(k+1)/2}\right)^{(k-2)/3k}\log^2 T\nonumber\\
&\ll y^2\log T+T^{(k+1)/3k}\log^2 T.\nonumber
\end{align}

Finally, from (5.1),(\ref{5.2}) and (5.9) we have
\begin{equation}
E_k(y,z)\ll y^2z^\varepsilon \log^4
T+T^{(k+1)/3k}z^\varepsilon\log^4 T.
\end{equation}

\section{\bf Proof of Theorem 1(Completion)}

First consider the case $k=4.$ Take $z=e^{10c_3\delta(T)},$ where
$c_3$ was the constant in (3.4). From (3.3) and (3.4) we get
\begin{eqnarray*}
\int_T^{2T}|R_2^{(4)}(x)+R_3^{(4)}(x)|^2dx\ll
T^{3/2}e^{-2c_3\delta(T)}.
\end{eqnarray*}
From (3.6)--(3.9), (4.4) and (5.10) we get
\begin{align*}
\int_T^{2T}|R_1^{(4)}(x)|^2dx=&\
\frac{B_4}{4\pi^2}\int_T^{2T}x^{1/2}dx
+O(T^{3/2}z_0^{-1/4}\log^{627} T)\\
&+O(Ty^2z^{\varepsilon}\log^5 T+T^{17/12}z^{\varepsilon}\log^6
T)\\
=&\
\frac{B_4}{4\pi^2}\int_T^{2T}x^{1/2}dx+O(T^{3/2}e^{-2c_3\delta(T)}).
\end{align*}

From the above two estimates and Cauchy's inequality we get
\begin{eqnarray*}
\int_T^{2T}R_1^{(4)}(x)(R_2^{(4)}(x)+R_3^{(4)}(x))dx\ll
T^{3/2}e^{-c_3\delta(T)}.
\end{eqnarray*}

From the above three estimates we get
\begin{align}
\int_T^{2T}|\Delta^{(4)}(x)|^2dx&=\int_T^{2T}|R_1^{(4)}(x)|^2dx+
2\int_T^{2T}R_1^{(4)}(x)(R_2^{(4)}(x)+R_3^{(4)}(x))dx\\
&+\int_T^{2T}|R_2^{(4)}(x)+R_3^{(4)}(x)|^2dx\nonumber\\
&=\frac{B_4}{4\pi^2}\int_T^{2T}x^{1/2}dx+O(T^{3/2}e^{-c_3\delta(T)}),\nonumber
\end{align}
which implies the case $k=4$ of Theorem 1.

Now suppose $k\geq 5.$ Take $z=T^{1-\varepsilon}.$ From (3.3) and
(3.5) we get
\begin{eqnarray*}
\int_T^{2T}|R_2^{(k)}(x)+R_3^{(k)}(x)|^2dx\ll
T^{1+\varepsilon}y^2+T^3y^{2-2k}.
\end{eqnarray*}
From (3.6)-(3.9), (4.4) and (5.10) we get
\begin{align*}
\int_T^{2T}|R_1^{(k)}(x)|^2dx =&\
\frac{B_k}{4\pi^2}\int_T^{2T}x^{1/2}dx
+O(T^{3/2+\varepsilon}y^{1/k-1/2})\\
&+O(T^{1+\varepsilon}y^2+T^{1+(k+1)/3k+\varepsilon}).
\end{align*}
The above two estimates imply
\begin{eqnarray*}
\int_T^{2T}R_1^{(k)}(x)(R_2^{(k)}(x)+R_3^{(k)}(x))dx\ll
T^{5/4+\varepsilon}y+T^{9/4}y^{1-k}.
\end{eqnarray*}

From the above three estimates we get
\begin{align*}
\int_T^{2T}|\Delta^{(k)}(x)|^2dx
=&\ \frac{B_k}{4\pi^2}\int_T^{2T}x^{1/2}dx+O(T^{1+(k+1)/3k+\varepsilon})\\
&+O(T^{5/4+\varepsilon}y+T^{9/4}y^{1-k}+T^{3/2+\varepsilon}y^{1/k-1/2}).
\end{align*}

Now on taking $y=T^{5/26}$ if $k=5$ and $y=T^{1/k-\varepsilon}$ if
$k\geq 6,$ we get
\begin{eqnarray}
\int_T^{2T}|\Delta^{(k)}(x)|^2dx
=\frac{B_k}{4\pi^2}\int_T^{2T}x^{1/2}dx+O(T^{\delta_k+\varepsilon}),
\end{eqnarray}
where $\delta_k$ was defined in Section 1. The case $k\geq 5$ of
Theorem 1  now follows from (6.2).

\section{\bf An expression of $\Delta(1,1,k;x)$ }
In order to prove Theorem 2, we shall give an expression of
$\Delta(1,1,k;x)$ in this section. We write
\begin{align}
D(1,1,k;x)=&\ \sum_{nd^k\leq x}d(n)\nonumber\\
=&\ \sum_{d\leq y}D(\frac{x}{d^k}) +\sum_{n\leq
\frac{x}{y^k}}d(n)\left[(\frac{x}{n})^{1/k}\right]-D\left(\frac{x}{y^k}
\right)[y]\nonumber\\
=&\ {\sum}_8+{\sum}_9-{\sum}_{10}\label{7.1}
\end{align}
say, where $x^{\varepsilon}\ll y\ll x^{1/k-\varepsilon}$ is a
parameter.

We write $\sum_8$ as
\begin{align*}
{\sum}_8=&\ \sum_{d\leq y}\left(\frac{x}{d^k}\log
\frac{x}{d^k}+(2\gamma-1)\frac{x}{d^k}+\Delta(\frac{x}{d^k})\right)\\
=&\ x\log x\sum_{d\leq y}\frac{1}{d^k}-kx \sum_{d\leq y}\frac{\log
d}{d^k}+(2\gamma-1)x\sum_{d\leq y}\frac{1}{d^k}+\sum_{d\leq
y}\Delta\left(\frac{x}{d^k}\right).
\end{align*}
By the well-known Euler-Maclaurin's formula we have
\begin{eqnarray*}
\sum_{d\leq y}\frac{1}{d^k}=\zeta(k)-\sum_{d>
y}\frac{1}{d^k}=\zeta(k)-\frac{y^{1-k}}{k-1}-\psi(y)y^{-k}+O(y^{-k-1})
\end{eqnarray*}
and
\begin{eqnarray*}
&&\sum_{d\leq y}\frac{\log d}{d^k}=-\zeta^{\prime}(k)-\sum_{d>
y}\frac{\log d}{d^k}\\&&=-\zeta^{\prime}(k)+\frac{y^{1-k}\log
y}{1-k}-\frac{y^{1-k}}{(k-1)^2}-\frac{\psi(y)\log
y}{y^{k}}+O(y^{-k-1}\log y).
\end{eqnarray*}

From the above three formulas we get
\begin{align}
{\sum}_8=&\ \zeta(k)x\log x-\frac{xy^{1-k}\log
x}{k-1}-\psi(y)xy^{-k}\log x\nonumber\\
&+k\zeta^{\prime}(k)x-\frac{kxy^{1-k}\log
y}{1-k}+\frac{kxy^{1-k}}{(k-1)^2}+\frac{kx\psi(y)\log
y}{y^{k}}\nonumber\\
&+(2\gamma-1)\zeta(k)x-(2\gamma-1)\frac{xy^{1-k}}{k-1}
-(2\gamma-1)\psi(y)xy^{-k}\nonumber\\
&+\sum_{d\leq y}\Delta\left(\frac{x}{d^k}\right)+O(xy^{-k-1}\log
x).
\end{align}

We write
\begin{align*}
{\sum}_9=&\ \sum_{n\leq
\frac{x}{y^k}}d(n)\left((x/n)^{1/k}-1/2-\psi((x/n)^{1/k})\right)\\
=&\ x^{1/k}\sum_{n\leq
\frac{x}{y^k}}d(n)n^{-1/k}-\frac{1}{2}D(xy^{-k})-\sum_{n\leq
\frac{x}{y^k}}d(n)\psi((x/n)^{1/k}).
\end{align*}
By partial summation we get($M=xy^{-k}$)
\begin{eqnarray*}
&&\sum_{n\leq M}d(n)n^{-1/k}=\int_{1^{-}}^M\frac{dD(u)}{u^{1/k}}=
\int_{1^{-}}^M\frac{d(u\log
u+(2\gamma-1)u)}{u^{1/k}}+\int_{1^{-}}^M\frac{d\Delta(u)}{u^{1/k}}\\
&&=\int_{1}^M\frac{\log
u+1+2\gamma-1}{u^{1/k}}du+\frac{\Delta(M)}{M^{1/k}}+\frac{1}{k}\int_{1}^M\frac{\Delta(u)}{u^{1+1/k}}du\\
&&=\zeta^2(1/k)+\frac{M^{1-1/k}\log
M}{1-1/k}-\frac{M^{1-1/k}}{(1-1/k)^2}+\frac{M^{1-1/k}}{1-1/k}+(2\gamma-1)\frac{M^{1-1/k}}{1-1/k}\\
&&\hspace{15mm}+\Delta(M)M^{-1/k}+O(M^{-1/k}),
\end{eqnarray*}
where we used the estimate
$$\int_M^\infty \frac{\Delta(u)}{u^{1+1/k}}du\ll M^{-1/k},$$
which follows from the well-known estimate $\int_1^t\Delta(u)du\ll
t.$

From the above two formulas we get
\begin{align}
{\sum}_9=&\ \zeta^2(1/k)x^{1/k}+\frac{xy^{1-k}\log
xy^{-k}}{1-1/k}-\frac{xy^{1-k}}{(1-1/k)^2}+\frac{xy^{1-k}}{1-1/k}\nonumber\\
&+(2\gamma-1)\frac{xy^{1-k}}{1-1/k}+y\Delta(xy^{-k})-\frac{1}{2}D(xy^{-k})\nonumber\\
&-\sum_{n\leq \frac{x}{y^k}}d(n)\psi((x/n)^{1/k})+O(y).
\end{align}

For $\sum_{10}$ we have
\begin{align}
-{\sum}_{10}=&\ \psi(y)xy^{-k}\log
xy^{-k}+(2\gamma-1)\psi(y)xy^{-k}+\psi(y)\Delta(xy^{-k})\\
&+\frac{1}{2}D(xy^{-k})-xy^{1-k}\log
xy^{-k}-(2\gamma-1)xy^{1-k}-y\Delta(xy^{-k}).\nonumber
\end{align}

From (7.1)--(7.4) we get
\begin{align*}
\Delta(1,1,k;x)=&\ \sum_{d\leq y}\Delta(\frac{x}{y^k})-\sum_{n\leq
\frac{x}{y^k}}d(n)\psi((x/n)^{1/k})+O(y)\\
&+O(xy^{-k-1}\log x)+O(|\Delta(xy^{-k})|).
\end{align*}

From $\Delta(u)\ll u^{1/3}$ we get
$$|\Delta(xy^{-k})|\ll x^{1/3}y^{-k/3}\ll y+xy^{-k-1}.$$

Thus we get the following Lemma .

{\bf Lemma 7.1.} Suppose $x^{\varepsilon}\ll y\ll
x^{1/k-\varepsilon}.$ Then
\begin{align*}
\Delta(1,1,k;x)=&\ \sum_{d\leq y}\Delta(\frac{x}{y^k})-\sum_{n\leq
\frac{x}{y^k}}d(n)\psi\left(\left(\frac{x}n\right)^{1/k}\right)
+O(xy^{-k-1}\log x)+O(y).
\end{align*}

\section{\bf Proof of Theorem 2}

It suffices for us to evaluate $\int_T^{2T}\Delta^2(1,1,k;x)dx$
for large $T.$ Suppose $T^{\varepsilon}\ll y\ll
T^{1/k-\varepsilon}$ is a parameter to be determined later and
$z=T^{1-\varepsilon}.$ For simplicity, we write ${\cal L}=\log T$
in this section. Similar to (3.1), by Lemma 7.1 we may write
\begin{equation}
\Delta(1,1,k;x)=R_{1,k}(x)+R_{2,k}(x)-R_{3,k}(x),
\end{equation}
where
\begin{eqnarray*}
&&R_{1,k}(x):=\frac{x^{1/4}}{\sqrt 2\pi}\sum_{d\leq
y}\frac{1}{d^{k/4}}\sum_{n\leq z}\frac{d(n)}{n^{3/4}}
\cos\left(4\pi\sqrt{\frac{nx}{d^k}}-\frac{\pi}{4}\right),\\
&&R_{2,k}(x):=\sum_{d\leq y}\Delta_2(\frac{x}{d^k};z),\\
&&R_{3,k}(x):=\sum_{n\leq\frac{x}{y^k}}d(n)\psi((x/n)^{1/k})
+O(xy^{-k-1}\log x)+O(y).
\end{eqnarray*}

Similar to the mean square of $R_1^{(k)}(x)$, we can prove that
\begin{align}
\int_T^{2T}|R_{1,k}(x)|^2dx =&\
\frac{C_k}{4\pi^2}\int_T^{2T}x^{1/2}dx
+O(T^{3/2+\varepsilon}y^{1/k-1/2})\nonumber\\
&+O(T^{1+\varepsilon}y^2+T^{1+(k+1)/3k+\varepsilon}).
\end{align}
From (3.3) we have
\begin{eqnarray}
\int_T^{2T}|R_{2k}(x)|^2dx\ll Ty^2{\cal L}^6.
\end{eqnarray}

Now we study the mean square of
$$S(x)=\sum_{n\leq\frac{x}{y^k}}d(n)\psi((x/n)^{1/k}).$$
Let $J=[\log^{-1} 2\log (Ty^{-k}{\cal L}^{-1})],$ then $J\ll {\cal
L}$ and we may write
\begin{eqnarray*}
&&S(x)=\sum_{j=0}^JS_j(x)+O({\cal L}^2),\\
&&S_j(x):=\sum_{x2^{-j-1}y^{-k}<n\leq
x2^{-j}y^{-k}}d(n)\psi((x/n)^{1/k}).
\end{eqnarray*}

Let $1/T\ll \eta<1/10$ is a real number and let $\eta T=N.$ Let
$$M(x,\eta):=\sum_{\eta x<n\leq 2\eta x}d(n)\psi((x/n)^{1/k}).$$
Then $S_j(x)=M(x,2^{-j-1}y^{-k}), j=0,1, \cdots, J$. We shall
study $\int_T^{2T}M^2(x,\eta)dx.$

According to Vaaler [16], we may write
$$\psi(t)=\sum_{1\leq |h|\leq N}a(h)e(ht)+O\left(\sum_{ |h|\leq
N}b(h)e(ht)\right)$$ with $a(h)\ll 1/|h|, b(h)\ll 1/N.$ Thus
\begin{eqnarray*}
M(x,\eta)&&=\sum_{1\leq |h|\leq N}a(h)\sum_{\eta x<n\leq 2\eta
x}d(n)e(h(x/n)^{1/k})\\
&&+O(\sum_{|h|\leq N}b(h)\sum_{\eta x<n\leq 2\eta
x}d(n)e(h(x/n)^{1/k}))\\
&&\ll 1+\sum_{1\leq h\leq N}h^{-1/2}h^{-1/2}\left|\sum_{\eta
x<n\leq 2\eta x}d(n)e(h(x/n)^{1/k})\right|.
\end{eqnarray*}
By Cauchy's inequality we get
\begin{eqnarray*}
M^2(x,\eta)\ll 1+\sum_{1\leq h\leq N}\frac{{\cal
L}}{h}\left|\sum_{\eta x<n\leq 2\eta
x}d(n)e(h(x/n)^{1/k})\right|^2.
\end{eqnarray*}
Integrating, squaring out and then by the first derivative test we
get
\begin{eqnarray*}
&&\int_T^{2T}M^2(x,\eta)dx\ll T+\sum_{1\leq h\leq N}\frac{{\cal
L}}{h}\int_T^{2T}\left|\sum_{\eta x<n\leq 2\eta
x}d(n)e(h(x/n)^{1/k})\right|^2dx\\
&&=T+\sum_{1\leq h\leq N}\frac{{\cal L}}{h}\int_T^{2T}\sum_{\eta
x<n\leq 2\eta x}d^2(n)dx\nonumber\\
&&\hspace{15mm}+\sum_{1\leq h\leq N}\frac{{\cal L}}{h}\int_T^{2T}
\sum_{\stackrel{\eta x<n,m\leq 2\eta x}{m\not= n}}
d(m)d(n)e(hx^{1/k}(m^{-1/k}-n^{-1/k}))dx\\
&&=O(TN{\cal L}^5)+\sum_{1\leq h\leq N}\frac{{\cal
L}}{h}\sum_{\stackrel{N<n,m\leq 4N}{m\not=
n}}d(m)d(n)\int_{I(m,n)}e(hx^{1/k}(m^{-1/k}-n^{-1/k}))dx\\
&&\ll TN{\cal L}^5+\sum_{1\leq h\leq N}\frac{{\cal
L}}{h}\sum_{\stackrel{N<n,m\leq 4N}{m\not=
n}}\frac{T^{1-1/k}d(n)d(m)}{h|m^{-1/k}-n^{-1/k}|}\\
&&\ll TN{\cal L}^5+\sum_{1\leq h\leq N}\frac{{\cal
L}}{h}\sum_{\stackrel{N<n,m\leq 4N}{m\not=
n}}\frac{T^{1-1/k}N^{1+1/k}d(n)d(m)}{h|m-n|}\\
&&\ll TN{\cal L}^5+T^{1-1/k}N^{2+1/k}{\cal L}^5,
\end{eqnarray*}
where $I(m,n)$ is a subinterval of $[T,2T].$

From Cauchy's inequality and the above estimate we get
\begin{eqnarray*}
&&\int_T^{2T}S^2(x)dx\ll \int_T^{2T}|\sum_{j=0}^JS_j(x)|^2dx+T{\cal L}^2\\
&&\ll {\cal L}\sum_{j=0}^J\int_T^{2T}|S_j(x)|^2dx+T{\cal L}^2\\
&&\ll (T^2y^{-k}+T^3y^{-2k-1}){\cal L}^6,
\end{eqnarray*}
which implies that
\begin{equation}
\int_T^{2T}R_{3k}^2(x)dx\ll (T^2y^{-k}+T^3y^{-2k-1}){\cal
L}^6+Ty^2.
\end{equation}

From (8.2)-(8.4) and Cauchy's inequality we get
\begin{equation}
\int_T^{2T}R_{1k}(x)(R_{2k}(x)+R_{3k}(x))dx\ll T^{5/4}y{\cal
L}^3+T^{7/4}y^{-k/2}{\cal L}^3+T^{9/4}y^{-k-1/2}{\cal L}^3.
\end{equation}
From (8.1)-(8.5) we get
\begin{eqnarray*}
&&\int_T^{2T}\Delta^2(1,1,k;x)dx=\frac{C_k}{4\pi^2}\int_T^{2T}
x^{1/2}dx\\&&\hspace{25mm}+O(T^{5/4}y{\cal
L}^3+T^{7/4}y^{-k/2}{\cal
L}^3+T^{9/4}y^{-k-1/2}{\cal L}^3)\\
&&\hspace{25mm}+O(T^{3/2}y^{1/k-1/2}{\cal
L}^{(k+1)^4+2}+T^{1+(k+1)/3k+\varepsilon}).
\end{eqnarray*}
Now on taking $y=T^{2/9}$ if $k=3,$ $y=T^{1/5}{\cal L}^{2496/5}$
if $k=4,$ $y=T^{5/26}{\cal L}^{10(6^4-1)/13}$ if $k=5$ and
$y=T^{1/k-\varepsilon}$ if $k\geq 6$ we get
\begin{align}
\int_T^{2T}\Delta^2(1,1,k;x)dx =&\ \frac{C_k}{4\pi^2}\int_T^{2T}
x^{1/2}dx +\left\{\begin{array}{ll}
O(T^{53/36}{\cal L}^3),&\mbox{if $k=3$, }\\
O(T^{29/20}{\cal L}^{503}),& \mbox{if $k=4,$}\\
O(T^{75/52}{\cal L}^{1000}),& \mbox{if $k=5,$}\\
O(T^{3/2-1/2k+1/k^2+\varepsilon}),& \mbox{if $k\geq 6.$}\\
\end{array}\right.
\end{align}
Theorem 2 follows from (8.6) immediately.


\begin{thebibliography}{99}

\bibitem{s1}R. C. Baker, The square-free divisor problem , Quart.
J. Math. Oxford Ser. (2) {\bf 45} (1994), no. 179, 269--277.

\bibitem{s2}R. C. Baker, The square-free divisor problem. II, Quart.
J. Math. Oxford Ser. (2) {\bf 47} (1996), no. 186, 133--146.

\bibitem{s3}E. Fouvry and H. Iwaniec, Exponential sums with
monomials, J. Number Theory {\bf 33} (1989), no. 3, 311--333.

\bibitem{s12}O. H\"older, \"Uber einen asymptotischen Ausdruck, Acta
Math. {\bf 59} (1932), 89--97.

\bibitem{s4}M. N. Huxley, Exponential sums and lattice points III,
Proc. London Math. Soc. {\bf87} (3) (2003), 591--609.

\bibitem{s5}A. Ivi\'c, The Riemann-zeta function, John Wiley \& Sons,
New York, 1985.

\bibitem{s6}A. Ivi\'c, The general divisor problem, J. Number
Theory {\bf 27} (1987), no. 1, 73--91.

\bibitem{s7}E. Kr\"atzel, Teilerprobleme in drei dimensionen,
Math. Nachr. {\bf 42} (1969), 275--288.

\bibitem{s8}A. Kumchev, The $k$-free divisor problem, Monatsh.
Math. {\bf 129} (2000), 321--327.

\bibitem{s9}F. Mertens, \"Uber einige asymptotische Gesetze der
Zahlentheorie, J. Reine. Angew. Math. {\bf 77} (1874), 289--338.

\bibitem{s10}T. Meurman, On the mean square of the Riemann zeta-function,
Quart. J. Math. Oxford Ser. (2) {\bf38} (1987), no. 151, 337--343.

\bibitem{s11}W. G. Nowak, M. Schmeier, Conditional asymptotic
formulae for a class of arithmetical functions, Proc. Amer. Math.
Soc. {\bf 103} (1988), 713--717.

\bibitem{s13}B. Saffari, Sur le nombre de diviseurs "$r$-libres" d'un
entier, et sur les points \`a coordonn\'ees enti\`eres dans
certaines r\'egions du plan, C. R. Acad. Sci. Paris. S\'er. A-B
{\bf 266} (1968), A601--A603.

\bibitem{s14}V. Siva Rama Prasad, D. Suryanarayana, The number of
$k$-free divisors of an integer, Acta Arith. {\bf 17} (1970/71),
345--354.

\bibitem{s15}K. C. Tong, On divisor problem III, Acta Math. Sinica {\bf 6}
(1956), 515-541.

\bibitem{s16}J. D. Vaaler, Some extremal functions in Fourier analysis, Bull.
Amer. Math. Soc. {\bf 12} (1985), 183--216.
\end{thebibliography}
\end{document}